\documentclass[12pt,a4paper,twoside]{article}
\usepackage{a4wide,amsmath,amsbsy,amsfonts,amssymb,stmaryrd}
\newcommand\Z{{\mathbb Z}}

\newcommand\Q{{\mathbb Q}}

\newcommand\C{{\mathbb C}}

\newcommand\cL{{\cal L}}

\newcommand\ra{\rightarrow}

\newcommand\GL{\mathrm{GL}}
\newcommand\SP{\mathrm{Sp}}

\newcommand\aut{\mathrm{Aut}}
\newcommand\out{\mathrm{Out}}

\newcommand\lra{\longrightarrow}
\newcommand\hookra{\hookrightarrow}
\newcommand\tura{\twoheadrightarrow}
\newcommand\da{\downarrow}

\newcommand{\pic}{\mathrm{Pic}}
 
\newcommand\sr{\stackrel}
\newcommand\st{\scriptstyle}
\newcommand\sst{\scriptscriptstyle}

\newcommand\ssm{\smallsetminus}
\newcommand\ol{\overline}
\newcommand\ccM{\overline{\cal M}}
\newcommand\cM{{\cal M}}
\newcommand\cA{{\cal A}}

\newcommand\cC{{\cal C}}

\newcommand\cG{{\cal G}}

\newcommand\wt{\widetilde}

\newcommand\wM{\wt{\cM}}

\newcommand\sg{\sigma}

\newtheorem{theorem}{Theorem}[section]
\newtheorem{corollary}[theorem]{Corollary}
\newtheorem{proposition}[theorem]{Proposition}
\newtheorem{lemma}[theorem]{Lemma}
\newtheorem{definition}[theorem]{Definition}
\newtheorem{remark}[theorem]{Remark}
\newenvironment{prf}[1]{\trivlist
\item[\hskip
\labelsep{\it #1.\hspace*{.3em}}]}{~\hspace{\fill}~$\square$\endtrivlist}
\newenvironment{proof}{\begin{prf}{Proof}}{\end{prf}}

\begin{document}

\title{Monodromy of stable curves of compact type: \\ rigidity and extension}
\author{Marco Boggi}\maketitle

\begin{abstract}Let $\cM_{g,n}$, for $2g-2+n>0$, be the moduli stack of $n$-pointed, genus $g$, smooth 
curves. For a family $C\ra S$ of such curves over a connected base and a geometric point
$\xi$ on $S$, the associated monodromy representation is the induced homomorphism
$\pi_1(S,\xi)\ra\pi_1(\cM_{g,n},[C_\xi])$ on algebraic fundamental groups. It is well known that, if $S$ is irreducible, 
reduced and locally of finite type over a field $k$ of characteristic zero, the fibre $C_\xi$
and the corresponding  monodromy representation determine the relative isomorphism class of the family.
In the first part of the paper, it is shown that suitable quotients of this representation suffice. 

These results are then applied to show that the monodromy representation associated to a family $C\ra S$
of $n$-pointed, genus $g$, stable curves of compact type, i.e. the induced homomorphism 
$\pi_1(S,\xi)\ra\pi_1(\wM_{g,n},[C_\xi])$ (where, $\wM_{g,n}$ denotes the moduli stack of $n$-pointed, genus $g$, 
stable curves of compact type), characterizes trivial and isotrivial families.

Let $U$ be an open subscheme of a normal, irreducible, locally noetherian scheme $S$ over a field $k$ of
characteristic zero and let $C\ra U$ be a family of stable curves of compact type. In the second part of
the paper, a monodromy criterion is given for extending $C\ra U$ to a family of stable curves of
compact type over $S$.    
\newline

\end{abstract}    



\section{Monodromy rigidity of families of curves}
According to a classical theorem by Grothendieck (see the theorem in the introduction of \cite{Gr}),
a family of abelian varieties in characteristic zero is determined by the associated monodromy
representation. More precisely, given two schemes (or stacks) $X$ and $X'$ over a common base $S$, we say that
they are {\it relatively isomorphic}, if they are isomorphic over $S$.
Then, it holds:

\begin{theorem}[Grothendieck]\label{mon gr}
Let $S$ be a reduced, connected scheme locally of finite type over a
field $k$ of characteristic zero, and let $\xi$ be a geometric point on $S$. The relative isomorphism class of a 
polarized abelian scheme $A$ over $S$, of relative dimension $g$, is determined by the fibre $A_\xi$ 
and by the corresponding monodromy representation $\pi_1(S,\xi)\ra\GL(H_1(A_\xi,\Z_\ell))$, for a 
prime $\ell>0$.
\end{theorem}

The corollary to the above theorem, given as well by Grothendieck in the introduction of \cite{Gr}, is better understood 
making use of the notion of fundamental group of an algebraic
stack (see \cite{N1}). As usual, the moduli stack of principally polarized abelian varieties of 
dimension $g$ is denoted by $\cA_g$.

\begin{corollary}[Grothendieck]\label{grothendieck} Let $S$ and $\xi$ be as above. A $k$-morphism $f:S\ra\cA_g$ is
determined by the image of the point $\xi$ and by the induced homomorphism
$f_\ast:\pi_1(S,\xi)\ra\pi_1(\cA_g, f(\xi))$.
\end{corollary}
\begin{proof}To the universal family ${\mathbb A}\ra\cA_g$ of abelian varieties and the point
$[A_\xi]\in\cA_g$ is associated the {\it universal monodromy representation}
$\pi_1(\cA_g,[A_\xi])\ra\GL(H_1(A_\xi,\Z_\ell))$. The monodromy representation
$\pi_1(S,\xi)\ra\GL(H_1(A_\xi,\Z_\ell))$ is the composition of $f_\ast$ with the universal monodromy
representation.
\end{proof}

\begin{remark}\label{oort}{\rm The above criterion implies that a family of abelian varieties over a simply
connected, reduced scheme $S$, over a field $k$ of characteristic zero, is trivial. Since this is false in
positive characteristic (see Remark~2.6 of \cite{Oort}), it follows, in particular, that the hypothesis on the
characteristic is essential in Theorem~\ref{mon gr} and most likely in all that follows in this section.}
\end{remark}

The Torelli morphism $t:\cM_g\rightarrow\cA_g$ is a natural morphism of D-M stacks fitting in a
commutative diagram:
$$\begin{array}{ccccc}
{\cal C}_\xi&\hookra&{\cal C}&\rightarrow&\cM_g\,\\
\,\,\da{\st j_{\sst P}}&&&&\,\,\da{\st t}\\
{\mathbb A}_{t(\xi)}&\hookra&{\mathbb A}&\ra&\cA_g,
\end{array}$$
where $\cC\ra\cM_g$ is the universal curve, $\cC_\xi$ is the fibre of $\cC$ over a geometric point
$\xi\in\cM_g$, ${\mathbb A}_{t(\xi)}$ is the fibre of $\cA_g$ over the point $t(\xi)\in\cA_g$ and
$j_P:\cC_\xi\hookra{\mathbb A}_{t(\xi)}$ is the embedding  of $\cC_\xi$ in its jacobian determined by a fixed point 
$P\in\cC_\xi$. Let us observe that, even though the embedding $j_p$ depends on $P$, the induced epimorphism
on fundamental groups $j_{P\ast}:\pi_1(\cC_\xi)\ra\pi_1({\mathbb A}_{t(\xi)})$ does not, since translations act trivially
on the geometric fundamental group of an abelian variety. Therefore, the abelianization
of the fundamental group of $\cC_\xi$ can be canonically identified with the fundamental group of its jacobian
${\mathbb A}_{t(\xi)}$. Since the geometric fundamental group of a smooth hyperbolic curve is center free, 
the universal curve $\cC\ra\cM_g$ determines on fundamental groups a short exact sequence (for more details, see 
Lemma~2.1 in \cite{M-T}):
$$1\ra\pi_1(\cC_\xi)\ra\pi_1(\cC)\ra\pi_1(\cM_g,\xi)\ra 1.$$
The associated universal monodromy representation then fits in the diagram: 
$$\begin{array}{cccc}
\pi_1(\cM_g,\xi)&\longrightarrow&\out (\pi_1(\cC_\xi))\,&\\
\,\,\da{\sst t_\ast}& &\da&(\ast)\\
\pi_1(\cA_g,t(\xi))&\longrightarrow&\GL(\pi_1({\mathbb A}_{t(\xi)})).&
\end{array}$$
Let us prove that the diagram $(\ast)$ is actually commutative. In
contrast with $\cC\ra\cM_g$, the universal curve $\cC'\ra\cM_{g,1}$ is endowed
with a tautological section $s:\cM_{g,1}\ra\cC'$. The Torelli morphism
$t':\cM_{g,1}\ra\cA_g$ and the Abel map $j_s:\cC'\ra{\mathbb A}$ associated to $s$ then fit in a natural commutative diagram:
$$\begin{array}{ccccc}
{\cal C'}_{\xi'}&\hookra&{\cal C'}&\rightarrow&\cM_{g,1}\,\\
\,\,\da{\st j_{\sst s(\xi')}}&&\,\,\da{\st j_{\sst s}}&&\,\,\da{\st t'}\\
{\mathbb A}_{t(\xi')}&\hookra&{\mathbb A}&\ra&\cA_g.
\end{array}$$
Therefore, the corresponding universal monodromy representations fit in a commutative diagram:
$$\begin{array}{ccc}
\pi_1(\cM_{g,1},\xi')&\longrightarrow&\aut (\pi_1(\cC'_{\xi'}))\,\\
\,\,\da{\sst t_\ast'}& &\da\\
\pi_1(\cA_g,t(\xi'))&\longrightarrow&\GL(\pi_1({\mathbb A}_{t(\xi')})).
\end{array}$$
Let us observe that the morphism $\cM_{g,1}\ra\cM_g$ (forgetting
the labeled point) is naturally identified with the universal curve $\cC\ra\cM_g$ and then $\cC'$ with
the pull-back of $\cC$ over $\cM_{g,1}$. Hence, taking the base point
$\xi'\in\cC_\xi$, the fibre $\cC'_{\xi'}$ is naturally isomorphic to
$\cC_\xi$ and there is a natural commutative diagram:
$$\begin{array}{ccc}
\pi_1(\cM_{g,1},\xi')&\longrightarrow&\aut (\pi_1(\cC'_{\xi'}))\,\\
\da& &\da\\
\pi_1(\cM_g,\xi)&\longrightarrow&\out (\pi_1(\cC_\xi)).
\end{array}$$
Since the homomorphism $t'_\ast$ factors through $t_\ast$, the diagram $(\ast)$
is commutative as well. From the above remarks, it now follows:

\begin{theorem}\label{homology rigidity}Let $S$ and $\xi\in S$ be a reduced, irreducible D-M stack locally of finite type over a
field $k$ of characteristic zero. Then, the relative isomorphism class of a smooth, proper, geometrically connected curve 
$C\ra S$ of genus $g\geq 2$ is determined by its fibre $C_\xi$ over $\xi$ and by the associated monodromy 
representation $\pi_1(S,\xi)\lra\GL(H_1(C_\xi,\Z_\ell))$, for a prime $\ell>0$.
\end{theorem}

\begin{proof}Let us consider first the case when $S$ is a normal irreducible scheme.
By Theorem~\ref{mon gr}, the monodromy representation $\pi_1(S,\xi)\lra\GL(H_1(C_\xi,\Z_\ell))$
determines the isomorphism class of the relative principally polarized Jacobian $J_{C/S}$. So,
it determines the isomorphism class of the principally polarized Jacobian $J_{C_{K(S)}}$ of the generic
fibre $C_{K(S)}$ of $C\ra S$ (where we denote by $K(S)$ the function field of $S$).

Torelli Theorem implies that the isomorphism class of a smooth hyperbolic proper curve over a
field of characteristic zero is determined by its Jacobian (for an elementary proof of this fact,
see Corollary~12.2 in \cite{Milne}). In particular, this is true for $C_{K(S)}$.

Since the D-M stack $\cM_g$ is separated and $S$ is normal, the isomorphism class of the curve $C\ra S$
is determined by the isomorphism class of its generic point. Hence,
the theorem follows for normal schemes. For reduced, irreducible D-M stacks, it will follow from:

\begin{lemma}\label{epimorphism}Let $\cM$ be a D-M stack endowed with a uniformisation $\cM'\ra\cM$, i.e. the stack $\cM'$ is
representable by a connected scheme and the given morphism is Galois {\'e}tale. 
Let then $f_i:X\ra \cM$, for $i=1,2$, be representable morphisms of D-M stacks, with $X$ irreducible, and $\pi:Y\ra X$ be a
representable epimorphism of D-M stacks. Then, $f_1\equiv f_2$ if and only if
$f_1\circ\pi=f_2\circ\pi$ and, for a geometric point $x\in X$, it holds 
$f_{1\ast}=f_{2\ast}:\pi_1(X,x)\ra\pi_1(\cM,f_1(x)=f_2(x))/N$, where $N$ is a normal subgroup containing the subgroup
of $\pi_1(\cM)$ corresponding to the cover $\cM'\ra\cM$ .
\end{lemma}
\begin{proof}One direction is obvious. So let us prove that, under the above hypotheses, $f_1=f_2$. 
By hypotheses, the two pull-backs of the Galois {\'e}tale cover $\cM'\ra\cM$ along $f_1$ and $f_2$ are naturally isomorphic.
Let us denote this {\'e}tale cover by $X'\ra X$ and let $f'_i:X'\ra\cM'$, for $i=1,2$, be the induced
morphisms of schemes. In order to prove that $f_1\equiv f_2$, it is enough to show that $f'_1=f'_2$.
Let $\pi':Y'\ra X'$ be the pull-back of the given epimorphism $\pi$ along $X'\ra X$, which is then an
epimorphism of schemes. By hypothesis, it holds $f'_1\circ\pi'=f'_2\circ\pi'$. Since $\pi'$ is an epimorphism
and all the spaces involved are schemes, it follows $f_1'=f_2'$ and hence $f_1\equiv f_2$.
\end{proof}

The moduli stack of curves $\cM_{g,n}$ is uniformisied by any abelian level structure $\cM^{(m)}_{g,n}$ with $m\geq 3$. 
In particular, since the given curve $C\ra S$ defines a representable morphism $S\ra\cM_{g,n}$, the D-M stack $S$ is
also uniformisable. Let then $S'$ be the normalization of a uniformisation of $S$. The natural morphism $S'\ra S$ is
clearly a representable epimorphism. So, we can apply Lemma~\ref{epimorphism}.   
Let then $C'\ra S'$ the pull-back of the given curve over $S'$. The homology monodromy representation associated to the curve 
$C'\ra S'$ is determined by that associated to $C\ra S$ and, as we proved above, the curve $C'\ra S'$
is determined by its associated monodromy representation. The theorem then follows. 
\end{proof}

\begin{remark}\label{warning}{\rm In the statement of Theorem~\ref{homology rigidity}, it is
essential to assume that $S$ is irreducible (in contrast with Theorem~\ref{mon gr}), 
as the following counter-
example (provided by the referee) shows.
Let $C\ra T$ be a smooth, projective, geometrically connected curve of genus $\geq 3$, over a smooth connected curve,
such that, generically, the group of automorphism is trivial and the fibre over $t\in T$ is hyperelliptic.
Let us glue $T$ with another copy of $T$ transversely in the points corresponding to $t$. Glueing the families above over $t$,
once with the identity and secondly with the hyperelliptic involution, we get two non-isomorphic families of curves.
However, the relative Jacobians are isomorphic. In fact, the hyperelliptic involution acts as $-1$ on the Jacobian
and this extends to all the relative Jacobian.}
\end{remark}

For a given algebraic curve $C$ and a prime $\ell\geq 2$, let us denote by $\pi_1^{(\ell)}(C)$ the pro-$\ell$ completion of its
fundamental group. Let us recall that the derived series of this group is defined inductively to be:
\[\begin{array}{lll}
\pi_1^{(\ell)}(C)(1)& := &\pi_1^{(\ell)}(C),\\
\pi_1^{(\ell)}(C)(k+1)& := &[\pi_1^{(\ell)}(C)(k),\pi_1^{(\ell)}(C)(k)].
\end{array}\]
Let us recall, moreover, that, for a topologically finitely generated pro-$\ell$ group, the terms of 
the derived series are closed subgroups.
  
Parshin's trick allows to extend Theorem~\ref{homology rigidity} to families of pointed curves:

\begin{theorem}\label{parshin}Let $S$ and $\xi\in S$ be as in the hypotheses of Theorem~\ref{homology rigidity}. 
Then, the relative isomorphism class of a smooth proper curve $C\ra S$ of genus $g$ with given 
non-coalescing sections $s_i: S\ra C$, for $i=1,\ldots, n$, such that $2g-2+n>0$,  is determined by its fibre $C_\xi$ 
over $\xi$ and by the induced monodromy representation:
$$\pi_1(S,\xi)\ra\out(\pi_1^{(\ell)}(C_\xi\ssm\cup_{i=1}^n
s_i(\xi))/\pi_1^{(\ell)}(C_\xi\ssm\cup_{i=1}^n s_i(\xi))(k)),$$
with $k=4$, if $n\geq 0$, and $k=3$, if $n\geq 2$.  
\end{theorem}

\begin{proof}Let us briefly recall the notion of $G$-Teichm{\"u}ller structure (for more details, see \cite{D-M}, \cite{P-dJ}
and \cite{B1}). Given an $n$-pointed, genus $g$, smooth curve $C\ra X$, with sections $s_i:X\ra C$, for $i=1,\ldots n$, 
a $G$-Teichm{\"u}ller structure is defined to be a global section of the sheaf whose stalk, at a geometric point 
$x\in X$, is the set of  exterior surjective group homomorphisms from $\pi_1(C_x\ssm\cup s_i(x),\tilde{x})$ 
to the group $G$. The stack of  $n$-pointed, genus $g$ smooth curves with $G$-Teichm{\"u}ller
structure is denoted by $\cM_{g,n}[G]$ and called the {\it D--M level structure} associated  to the group $G$.
It is endowed with a natural Galois {\'e}tale morphism $\cM_{g,n}[G]\ra\cM_{g,n}$

For any finitely generated group $G$, the subgroup $G^2$ generated by $2$-nd order powers of $G$ is a finite index 
verbal subgroup (in particular characteristic)
containing the second term of the derived series $G(2)$. So, its commutator subgroup $[G^2,G^2]$ contains $G(3)$. 
Going one step further, $(G^2)^2$ is a finite index verbal subgroup of $G$ containing $G(3)$. 
Hence its commutator subgroup $[(G^2)^2,(G^2)^2]$ contains $G(4)$.
  
In order to prove the theorem, we need a characteristic subgroup $K$ of $\pi_1^{(\ell)}(C_\xi\ssm\cup_{i=1}^n s_i(\xi))$, such that 
the corresponding Galois cover $C^K\ra C_\xi$ branches non-trivially over the points $s_i(\xi)$, for $i=1,\ldots, n$, and the genus of 
$C^K$ is at least two. For $n\geq 2$, the subgroup $\pi_1^{(\ell)}(C_\xi\ssm\cup_{i=1}^n s_i(\xi))^2$ would already do but, for 
$n=1$, we need to take $K:=(\pi_1^{(\ell)}(C_\xi\ssm\cup_{i=1}^n s_i(\xi))^2)^2$.

Let then $G:=\pi_1^{(\ell)}(C_\xi\ssm\cup_{i=1}^n s_i(\xi))/K$. By the theory
of $G$-Teichm{\"u}ller structures, it follows that there is a finite connected, Galois {\'e}tale cover $S'\ra S$ such that
the pull-back $C'\ra S'$ of the given curve is an $n$-pointed curve, with sections $s'_i$, for $i=1,\ldots, n$, and is 
endowed with a $G$-Teichm{\"u}ller structure, such that the induced
relative Galois cover $\cC^K\ra C'$ over $S'$ has fibrewise the required properties.

Let $\xi'\in S'$ be a geometric point lying over $\xi$. The fibre $\cC^K_{\xi'}$ is naturally isomorphic to $C^K$.
By Theorem~\ref{homology rigidity}, the curve $\cC^K\ra S'$ is determined by the 
representation $\pi_1(S',\xi')\ra\GL(H_1(C^K,\Z_\ell))$ and the fibre $C^K$. In their turn, they are determined by
the monodromy representation given in the statement of the theorem and by the fibre $C_\xi$.

Now, the curve $\cC^K\ra S'$ determines its quotient 
$C'\ra S'$ and the set of sections $\{s_i'\}$ as the connected components of the branch locus of the natural morphism 
$\cC^K\ra C'$. The order on the latter set is then determined by the order on the points $s(\xi_i)$ on $C_\xi$ given by the
indices. At this point, to obtain the theorem, we just apply Lemma~\ref{epimorphism}.

\end{proof}


If one is just interested in knowing whether the given $n$-pointed curve $C\ra S$ is isotrivial (or just trivial),
Theorem~\ref{parshin} can be considerably improved. Let $C_{g,n}$ be an $n$-punctured algebraic curve of genus $g$ 
over a field $k$ of characteristic zero and let $C_g:=C_{g,0}$ be its compactification. Let us consider the following natural
filtrations on its algebraic fundamental group $\pi_1(C_{g,n})$. The descending central series is defined by
$\pi_1(C_{g,n})^{[1]}:=\pi_1(C_{g,n})$ and $\pi_1(C_{g,n})^{[k]}:=[\pi_1(C_{g,n})^{[k-1]},\pi_1(C_{g,n})]$. The
weight filtration, which comes from Hodge theory, is defined as follows. Let $N$ be the kernel of the 
natural morphism $\pi_1(C_{g,n})\ra\pi_1(C_{g})$ (filling in the punctures) and define
\[\begin{array}{lll}
W^1 \pi_1(C_{g,n})& := &\pi_1(C_{g,n}),\\
W^2 \pi_1(C_{g,n})& := &N\cdot \pi_1(C_{g,n})^{[2]},\\
W^{k+1} \pi_1(C_{g,n})& := &[\pi_1(C_{g,n}),W^{k} \pi_1(C_{g,n})]\cdot [N, W^{k-1}\pi_1(C_{g,n})],
\end{array}\]
where we denote by $\pi_1(C_{g,n})^{[k]}$ the $k$-th term of the descending central series, i.e. the
normal closed subgroup spanned by commutators of order $k$. 
As for the descending central series, it holds $[W^s \pi_1(C_{g,n}),W^t \pi_1(C_{g,n})]\leq W^{s+t}
\pi_1(C_{g,n})$. The descending central series and the weight filtration are cofinal to each other (for $n=0$,
they coincide). In fact, it holds $W^{2k-1}\pi_1(C_{g,n})\leq\pi_1(C_{g,n})^{[k]}\leq W^k \pi_1(C_{g,n})$. Similar
filtrations can then be defined for the pro-$\ell$ completion $\pi_1^{(\ell)}(C_{g,n})$, for a prime $\ell\geq 0$.

By the results in \cite{B}, the quotient of the fundamental group of $C_{g,n}$ by the third term of the weight
filtration has a remarkable geometric meaning in the theory of moduli of curves. 
Let $\cC\ra\cM_{g,n}$ be the $n$-punctured, genus $g$ curve obtained from the universal curve removing the images
of the tautological sections (i.e. the curve canonically isomorphic to $\cM_{g,n+1}$)  
and let $\wM_{g,n}$ be the moduli stack of $n$-pointed, genus
$g$ stable curves of compact type. According to the results of \S 3 in \cite{B}, the universal
monodromy representation $\pi_1(\cM_{g,n},\xi)\ra\out(\pi_1(\cC_\xi))$
induces a representation:
$$\pi_1(\wM_{g,n},\xi)\ra\out(\pi_1(\cC_\xi)/W^3\pi_1(\cC_\xi)),$$
whose restriction to the geometric fundamental group is faithful. We can further specialize this
representation to a given prime $\ell$. As remarked in \S 3 of \cite{B}, the profinite group $\pi_1(C_\xi)/W^3\pi_1(C_\xi)$ 
decomposes as the direct product of its maximal pro-$\ell$ quotients. Therefore, there is a natural isomorphism:
$$\out(\pi_1(C_\xi)/W^3\pi_1(C_\xi))\cong\prod_{\ell\;\mathrm{prime}}\out(\pi_1^{(\ell)}(C_\xi)/W^3\pi_1^{(\ell)}(C_\xi)).$$
Let us then denote by $\pi_1^{(\ell)}(\wM_{g,n},\xi)$, for a given prime $\ell>0$ the image of $\pi_1(\wM_{g,n},\xi)$,
in the virtual pro-$\ell$ group $\out(\pi_1^{(\ell)}(C_\xi)/W^3\pi_1^{(\ell)}(C_\xi))$, by the monodromy representation.

The (iso)-triviality of an $n$-pointed, genus $g$ relative curve is controlled by the monodromy
representation with values in this group. More precisely, it holds:

\begin{proposition}\label{parshin improved}Let $S$ be a connected, reduced D-M stack, locally of finite 
type over a field $k$ of characteristic zero, and $\xi\in S$ a geometric point. 
For $2g-2+n>0$ and $g\geq 1$, an $n$-punctured smooth curve $C\ra S$ of genus $g$ is isotrivial (respectively trivial)
if and only if the associated monodromy representation $\pi_1(S,\xi)\ra\out(\pi_1^{(\ell)}(C_\xi)/W^3\pi_1^{(\ell)}(C_\xi))$
has finite (respectively trivial) image, for a prime $\ell>0$. 
\end{proposition}
\begin{proof}The only if statements are obvious. We can clearly assume that $S$ is irreducible. 
Moreover, possibly replacing $S$ by a Galois, {\'e}tale cover $S'\ra S$,
it is clealy enough to prove just the triviality statement for $S$ an irreducible, reduced scheme. 

The case $g\geq 2$ and $n=0$ is dealt by Theorem~\ref{homology rigidity} while the case $g=1$ and $n=1$ is dealt by 
Theorem~\ref{mon gr}. Then, the proof proceeds by induction on $n$. So, let us assume that the proposition has been 
proved for $n$ punctures and let us prove the $(n+1)$-punctured case.  

Let then $C\ra S$ be an $(n+1)$-punctured curve of genus $g$ over a reduced, irreducible scheme and let $C'\ra S$
be the $n$-punctured curve obtained from $C\ra S$ filling in the $(n+1)$-th section. The natural embedding
$C_\xi\hookra C_\xi'$ induces an epimorphism on truncated fundamental
groups so that the monodromy representation $\pi_1(S,\xi)\ra\out(\pi_1^{(\ell)}(C'_\xi)/W^3\pi_1^{(\ell)}(C'_\xi))$ is
induced by the monodromy representation $\pi_1(S,\xi)\ra\out(\pi_1^{(\ell)}(C_\xi)/W^3\pi_1^{(\ell)}(C_\xi))$. 
Therefore, since the latter is trivial, the former is trivial as well. So, by inductive hypothesis, the family
of $n$-punctured curves $C'\ra S$ is trivial.  

Let $f:S\ra\wM_{g,n+1}$ be the morphism corresponding to the given curve $C\ra S$. The image of $f$ is then contained
in the fibre of the curve $\wM_{g,n+1}\ra\wM_{g,n}$ over the geometric point which parametrizes the curve $C'_\xi$. Let us denote this point by $\zeta$ and the fibre above it by $F_\zeta$. By Theorem~3.11 in \cite{B}, there is a short exact sequence:
$$1\ra H_1(F_\zeta,\Z_\ell)\ra\pi_1^{(\ell)}(\wM_{g,n+1},f(\zeta))\ra\pi_1^{(\ell)}(\wM_{g,n},\zeta)\ra 1.$$
As we already remarked, there is a natural representation 
$$\pi_1^{(\ell)}(\wM_{g,n+1},f(\xi))\ra\out(\pi_1^{(\ell)}(C_\xi)/W^3\pi_1^{(\ell)}(C_\xi))$$ 
whose restriction to the geometric fundamental group is faithful. In particular, by hypothesis, the induced homomorphism 
$f_\ast: H_1(S\otimes\ol{k},\Z_\ell)\ra H_1(F_\zeta,\Z_\ell)$ is then trivial. Since $F_\zeta$ is a smooth curve of genus 
$\geq 1$, it follows that $f$ is the constant map.
\end{proof}


The above description of the fundamental group of the moduli stack $\wM_{g,n}$ suggests that a criterion similar to the 
one given in  Proposition~\ref{parshin improved} could still hold for families of curves of compact type. 
Let us remark that the open substack $\wM_{g,n}$ of $\ccM_{g,n}$ is the largest for which a monodromy criterion to 
characterize morphisms $S\ra\wM_{g,n}$ makes sense. In fact, by Teichm{\"u}ller theory, the topological fundamental group 
of the stack $\cM_{g,n}$, for $g\geq 1$, is normally generated by any ``small'' loop around the divisor whose generic point 
parameterizes singular irreducible curves. Therefore, an open substack of $\ccM_{g,n}$, containing a point parameterizing a 
singular irreducible curve, is simply connected. 

However, it is clear that a space, which enjoys the same
kind of rigidity property stated above, cannot contain any rational curve, while this not the case, in general,
for $\wM_{g,n}$ (this is why $\wM_{g,n}$ fails to be,  also topologically, an Eilenberg-Maclane
space of type $(\pi,1)$). So, we will obtain a criterion actually weaker than Proposition~\ref{parshin improved}.

Let us then give the following definition:
\begin{definition}\label{linearly rigid}{\rm Let $2g-2+n>0$. A $n$-pointed, genus $g$ stable
curve of compact type $C\ra S$, with labeling sections $s_i:S\ra C$ for $i=1,\ldots n$,  is {\it linearly
rigid} if, for every fibre $C_\xi$ over a geometric point $\xi\in S$, the following conditions are satisfied: 
\begin{enumerate}
\item the number of points where an irreducible rational component $E$ of $C_\xi$ meets the other
components of the curve $C_\xi$ and the images of the sections $s_i$, for $i=1,\ldots n$, is exactly $3$; 
\item if $E_1$ and $E_2$ are distinct irreducible rational components of $C_\xi$, then $E_1\cap E_2=\emptyset$.
\end{enumerate}
Instead, we will say that the curve $C\ra S$ is {\it linear} if its moduli only depend on the rational components of the fibers}
\end{definition}

The clutching morphisms $\wM_{g_1,n_1}\times\ldots\times\wM_{g_k,n_k}\ra\wM_{g,n}$, obtained by restriction from 
those on $\ccM_{g,n}$, are proper. Therefore, their images are closed substacks of 
$\wM_{g,n}$. The stack of linearly rigid curves is precisely the open substack of $\wM_{g,n}$, 
whose complement is the union of all closed strata such that the 
domain of the corresponding clutching morphisms contains a factor of the type $\wM_{0,k_i}$, with $k_i>3$. 
It does not contain any obvious rational curve. We will see that it does not contain rational curves at all.

\begin{theorem}\label{triviality}Let $S$ be a connected, reduced D-M stack, locally of finite 
type over a field $k$ of characteristic zero, and $\xi$ be a geometric point on $S$. Let then $C\ra S$, for $2g-2+n>0$, be an
$n$-pointed, genus $g$, stable curve of compact type. If the monodromy representation
$\pi_1(S,\xi)\ra\pi_1^{(\ell)}(\wM_{g,n},\xi)$ is finite (respectively trivial), for a prime $\ell>0$, the curve $C\ra S$ 
is a linear family. In particular, if, for some geometric point $\xi\in S$, the fibre $C_\xi$ is 
linearly rigid, the curve $C\ra S$ is an isotrivial family (respectively a trivial family).
\end{theorem}

\begin{proof}It is clearly enough to prove the theorem for the case when the monodromy representation is trivial and $S$ is an
irreducible scheme. Let $f:S\ra\wM_{g,n}$ be the morphism induced by $C\ra S$. If $f$ maps the generic point of $S$
to $\cM_{g,n}$, the conclusion follows from Proposition~\ref{parshin improved}. Otherwise, possibly
after replacing $S$ by an {\'e}tale cover, we may assume that $f$ factors through a certain clutching
morphism $\prod_{i=1}^k\wM_{g_i,n_i+1}\ra\wM_{g,n}$, with $\sum_{i=1}^k g_i=g$ and $\sum_{i=1}^k
n_i=n$. Let the factoring morphism $f':S\ra\prod_{i=1}^k\wM_{g_i,n_i+1}$ be optimal, in the sense that,
composed with any of the projections $p_i:\prod_{i=1}^k\wM_{g_i,n_i+1}\ra\wM_{g_i,n_i+1}$, it maps the 
generic point of $S$ inside $\cM_{g_i,n_i+1}$, for $i=1,\ldots k$.

If $g_i\neq 0$, by Proposition~\ref{parshin improved}, the morphism $p_i\circ f'$ is constant whenever the monodromy 
representation $(p_i\circ f')_\ast:\pi_1(S,\xi)\ra\pi_1^{(\ell)}(\wM_{g_i,n_i+1},p_i\circ f'(\xi))$ is trivial.
On the other hand, by Theorem~4.2 in \cite{B}, we have that the homomorphism
$(p_i\circ f')_\ast$ is determined by the monodromy representation
$f_\ast:\pi_1(S,\xi)\ra\pi_1^{(\ell)}(\wM_{g,n},[C_\xi])$ which, by hypothesis is trivial. Therefore, $(p_i\circ f')_\ast$ is trivial as well.
In conclusion, the morphism $p_i\circ f'$ is constant whenever $g_i\neq 0$. This completes the proof of the theorem.
\end{proof}

\begin{corollary}\label{simply connected}Let $S$ and $\xi\in S$ be as in the hypotheses
of Theorem~\ref{triviality} and let $C\ra S$, for $2g-2+n>0$, be an $n$-pointed, genus $g$, stable
curve of compact type. If $S$ is simply connected, the curve $C\ra S$ is a linear family. In particular, if, for
some geometric point $\xi\in S$, the fibre $C_\xi$ is linearly rigid, the curve $C\ra S$ is a trivial family.
\end{corollary}

\section{Extending families of curves}\label{extension}
Theorem~\ref{homology rigidity} and Theorem~\ref{parshin improved} show to which extent the monodromy
representation characterizes a family of smooth curves. Suppose given a dense open subscheme $U$ of a
scheme $S$ as above and an $n$-punctured smooth curve $C\ra U$ and that the monodromy representation
$\pi_1(U,\xi)\ra\out(\pi_1(C_\xi))$ "extends" to a representation  $\pi_1(S,\xi)\ra\out(\pi_1(C_\xi))$. Is it
then possible to extend $C\ra U$ to an $n$-punctured smooth curve $C\ra S$?

In \cite{Stix1}, Stix answered positively to this question with some restriction on $S$ but with no
restriction on the characteristic:

\begin{theorem}[Stix]\label{stix}Let $2g-2+n>0$. Let $U$ be a dense open subscheme of a normal,
connected, excellent, locally noetherian scheme $S$, $\xi$ a geometric point of $U$ and $\ell$ a prime
number invertible on $S$. Let $C\ra U$ be an $n$-punctured smooth curve of genus $g$. The curve $C\ra U$
extends to an $n$-punctured smooth curve of genus $g$ over $S$, if and only if the monodromy
representation $\pi_1(U,\xi)\ra\out(\pi_1^{(\ell)}(C_\xi))$ factors through the epimorphisms
$\pi_1(U,\xi)\tura\pi_1(S,\xi)$, induced by the injection $U\subset S$, and a homomorphism
$\pi_1(S,\xi)\ra\out(\pi_1^{(\ell)}(C_\xi))$.
\end{theorem}

Here, we will prove the following refinement of Stix' Theorem:

\begin{theorem}\label{ext}Let $2g-2+n>0$. Let $U\subset S$, $\xi$ and $\ell$ as in the hypotheses of
Theorem~\ref{stix}. Let $C\ra U$ be an $n$-punctured smooth curve of genus $g$. The curve $C\ra U$ 
extends to an $n$-punctured smooth curve of genus $g$ over $S$, if and only if the monodromy
representation $\pi_1(U,\xi)\ra\out(\pi_1^{(\ell)}(C_\xi)/W^{4}\pi_1^{(\ell)}(C_\xi))$ factors through the
epimorphisms $\pi_1(U,\xi)\tura\pi_1(S,\xi)$, induced by the inclusion $U\subset S$, and a homomorphism
$\pi_1(S,\xi)\ra\out(\pi_1^{(\ell)}(C_\xi)/W^{4}\pi_1^{(\ell)}(C_\xi))$.
\end{theorem}

\begin{proof}The proof (with the exception of Lemma~\ref{cohomology}) closely follows \cite{Stix1} and
\cite{Stix2}. As in \cite{Stix1}, the first step is the following:

\begin{lemma}\label{1}Let $U$, $S$, $\xi$ and $\ell$ be as in the statement of Theorem~\ref{ext}. Let
$f:U\ra\cM_{g,n}$, for $2g-2+n>0$, be a map such that the associated monodromy representation
$\pi_1(U,\xi)\ra\out(\pi_1^{(\ell)}(C_\xi)/W^{4}\pi_1^{(\ell)}(C_\xi))$ factors through the epimorphisms
$\pi_1(U,\xi)\tura\pi_1(S,\xi)$, induced by the inclusion $U\subset S$. Then, there exists a projective
birational morphism $\sg:S'\ra S$ which is an isomorphism above $U$, such that:
\begin{enumerate}
\item $S'$ is normal;
\item $\sg_\ast {\cal O}_{S'}={\cal O}_S$;
\item the map $f$ extends to a map $f':S'\ra\cM_{g,n}$, such that the associated monodromy
representation $\pi_1(S',\xi)\ra\out(\pi_1^{(\ell)}(C_\xi)/W^{4}\pi_1^{(\ell)}(C_\xi))$ factors as:
$$\pi_1(S',\xi)\sr{\sg_\ast}{\ra}\pi_1(S,\xi)\sr{f_\ast}{\ra}
\out(\pi_1^{(\ell)}(C_\xi)/W^{4}\pi_1^{(\ell)}(C_\xi)).$$
\end{enumerate}
\end{lemma}

\begin{proof}Let $\cM^{(\ell)}\ra\cM_{g,n}$ be the abelian level structure of order $\ell$ (for $\ell=2$, 
consider instead the abelian level $\cM^{(4)}\ra\cM_{g,n}$ and change the rest accordingly). It is a
quasi-projective variety, representable over $\cM_{g,n}$. Let then
$S(\ell):=S\times_{\cM_{g,n}}\cM^{(\ell)}$. The natural map
$f(\ell): S(\ell)\ra\cM^{(\ell)}$ is $\SP_{2g}(\Z/\ell)$-equivariant. Since the action of $\SP_{2g}(\Z/\ell)$ on
$\cM^{(\ell)}$ extends to its Deligne-Mumford compactification $\ccM^{(\ell)}$, the natural action of
$\SP_{2g}(\Z/\ell)$ on the graph $\cG_{f(\ell)}$ of $f(\ell)$ extends to its closure $\ol{\cG_{f(\ell)}}\subset
S(\ell)\times\ccM^{(\ell)}$. Let then $S'$ be the normalization of the quotient
$\ol{\cG_{f(\ell)}}/\SP_{2g}(\Z/\ell)$. The natural morphism $\sg:S'\ra S$ has all the required properties. The
only non-trivial property to check is that the image of the natural morphism $f':S'\ra\ccM_{g,n}$ is actually
contained in $\cM_{g,n}$. Here, the argument goes exactly as in the proof of Lemma~2.3 in \cite{Stix2}, 
since the image of $\pi_1(\cM_{g,n})$ in $\out(\pi_1^{(\ell)}(C_\xi)/W^{4}\pi_1^{(\ell)}(C_\xi))$ classifies a
tower of Galois covers of $\ccM_{g,n}$ which ramifies universally along the Deligne-Mumford boundary 
(see \cite{B-P} and \cite{Pikaart}).
\end{proof} 
 
The difficult part of the proof of Theorem~\ref{ext} consists in showing that the morphism constructed
above $f':S'\ra\cM_{g,n}$ actually factors through $\sg:S'\ra S$. Thanks to property $ii.)$ of Lemma~\ref{1},
it is enough to prove that this is true set-theoretically, i.e. that the fibers of the morphism $\sg$ are
contracted by $f'$. By construction, the fibers of $\sg$ are proper sub-schemes of $S'$. Therefore, it is now
enough to prove:

\begin{lemma}\label{2}Let $C\ra S$, for $2g-2+n>0$, be an $n$-punctured smooth curve of genus $g$ over a
proper smooth connected curve over an algebraically closed field $\ol{k}$, 
such that the associated monodromy representation
$\rho_S:\pi_1(S,\xi)\ra\out(\pi_1^{(\ell)}(C_\xi)/W^{4}\pi_1^{(\ell)}(C_\xi))$ is trivial. 
Then, the curve $C\ra S$ is a trivial family.
\end{lemma}

\begin{proof}In the previous section, we actually proved a stronger result but with the further assumption
that $C$ be defined over a field of characteristic zero. So, we need to provide a different proof.

For $g\leq 1$, the lemma follows trivially from the fact that the coarse moduli spaces
$M_{g,n}$, in these cases, are affine varieties and do not contain proper curves. For the
case $g\geq 2$, we follow the strategy of Proposition~2.5 in \cite{Stix2}. 

Let $\cL$ be Mumford's ample line bundle on $\cM_{g,n}\otimes\ol{k}$ and let
$k_1\in H^2(\cM_{g,n}\otimes\ol{k},\Q_\ell)$ be its first Chern class. If $f:S\ra\cM_{g,n}$ is not the constant map, the
pull-back $f^\ast(\cL)$ of $\cL$ over $S$ is ample as well and, in particular, one must have
$c_1(f^\ast(\cL))\neq 0$. There are natural representations
$\rho_k:\pi_1(\cM_{g,n})\ra\out(\pi_1^{(\ell)}(C_\xi)/W^{k}\pi_1^{(\ell)}(C_\xi))$, for $k\geq 1$, induced by
the universal monodromy representation. Let us denote by $G_{k,\ell}$ the image of $\pi_1(\cM_{g,n}\otimes\ol{k})$ in 
the profinite group $\out(\pi_1^{(\ell)}(C_\xi)/W^{k}\pi_1^{(\ell)}(C_\xi))$. Consider then the
commutative diagram:
$$
\begin{array}{ccccc}
\pic(\cM_{g,n}\otimes\ol{k})&\sr{c_1}{\ra}& H^2(\cM_{g,n}\otimes\ol{k},\Q_\ell)&\sr{\rho_4^\ast}{\leftarrow}& H^2(G_{4,\ell},\Q_\ell)\,\\
\da{\st f^\ast}&&\da{\st f^\ast}&&\da{\st\rho_S^\ast}\\
\pic(S)&\sr{c_1}{\ra}& H^2(S,\Q_\ell)&\sr{\sim}{\leftarrow}& H^2(\pi_1(S),\Q_\ell).
\end{array}
$$
Lemma~\ref{2} now follows from:

\begin{lemma}\label{cohomology} For $g\geq 2$, the image of $\rho_4^\ast$ in $H^2(\cM_{g,n}\otimes\ol{k},\Q_\ell)$
contains the Mumford class $k_1$.
\end{lemma}
\begin{proof}Let us consider the following commutative diagram:
$$\begin{array}{ccc}
\pi_1^{\mathrm{top}}(\cM_{g,n}\otimes\C)&\sr{\rho_4}{\ra} &G_{4,\ell}\,\\
&\searrow^{\rho_3}&\da\\
&&G_{3,\ell}.
\end{array}$$
By Proposition~5.1 in \cite{B}, the homomorphism
$$\rho_3^\ast: H^2(G_{3,\ell},\Q_\ell)\ra H^2(\pi_1^{\mathrm{top}}(\cM_{g,n}\otimes\C),\Q_\ell)\cong
H^2(\cM_{g,n}\otimes\C,\Q_\ell)$$ 
is surjective. Therefore, by the above diagram, the same holds for $\rho_4^\ast$.
\end{proof}

In fact, from the hypothesis that $\rho_S$ is trivial, it follows that $\rho_S^*$ is the zero map. Therefore,
$c_1(f^*(\cL))=0$ and then $f$ is necessarily the constant map.
\end{proof}
\end{proof}

In the light of Theorem~\ref{triviality}, it is natural to expect that a result of the type of
Theorem~\ref{ext} holds for families of stable curves of compact type as well. In this section, we will show
that this is actually the case, at least after some restrictions on the data. The starting point is the following
theorem by Grothendieck (see \S 4 in \cite{Gr}):

\begin{theorem}[Grothendieck]\label{ext grothendieck}Let $S$ be a normal, irreducible, locally noetherian scheme 
over a field $k$ of characteristic zero, and let $U$ be an open subscheme of $S$ and $\xi\in U$ a geometric point.
Let $f:U\ra\cA_g$ be a $k$-morphism. Then, $f$ extends to a $k$-morphism $S\ra\cA_g$, if and only if the
induced homomorphism $f_\ast:\pi_1(U,\xi)\ra\pi_1(\cA_g,f(\xi))$ factors through the
epimorphism $\pi_1(U,\xi)\tura\pi_1(S,\xi)$, induced by the injection $U\hookra S$, and a homomorphism
$\pi_1(S,\xi)\ra\pi_1(\cA_g,f(\xi))$.
\end{theorem}

Let us suppose given schemes $U\subset S$ as above and a $k$-morphism  $f:U\ra\wM_{g,n}$ such 
that the induced homomorphism $f_\ast:\pi_1(U,\xi)\ra\pi_1(\wM_{g,n},f(\xi))$ factors
through the epimorphism $\pi_1(U,\xi)\tura\pi_1(S,\xi)$ induced by the injection $U\subset S$. There is a
natural $k$-morphism $t:\wM_{g,n}\ra\cA_g$. The composition $t\circ f: U\ra\cA_g$ satisfies the
conditions of Theorem~\ref{ext grothendieck} and hence extends to a $k$-morphism $S\ra\cA_g$.

Let $C\subset S$ be a curve which intersects non-trivially $U$. The image of the extension $C\ra\ccM_{g,n}$ of the
induced morphism $f:C\cap U\ra\wM_{g,n}$ must then be contained in $\wM_{g,n}$, i.e. $f:C\cap U\ra\wM_{g,n}$
extends to a morphism $\ol{f}:C\ra\wM_{g,n}$.

For any given point $P\in S\ssm U$, there is a curve $C\subset S$ as above such that $P\in C$.
{\it The set of limit curves of $U$ in $P$} is then defined to be the set of
isomorphism classes of $n$-pointed, genus $g$, stable curves of compact type parametrized by the 
images of $P$ by all the possible morphisms $C\ra\wM_{g,n}$ as above. Let us remark that the jacobians  
of the limit curves of $U$ in $P$ are all isomorphic. 

We can now state the monodromy extension criterion for families of stable curves of compact type:  

\begin{theorem}\label{ext compact type}Let $S$ be a normal, irreducible, locally noetherian scheme over a 
field $k$ of characteristic zero, and let $U$ be an open subscheme of $S$ and $\xi\in U$ a geometric point.
Let $f:U\ra\wM_{g,n}$, for $2g-2+n>0$, be a $k$-morphism. Let us assume that the induced homomorphism
$f_\ast:\pi_1(U,\xi)\ra\pi_1(\wM_{g,n},f(\xi))$ factors through the epimorphism
$\pi_1(U,\xi)\tura\pi_1(S,\xi)$, induced by the injection $U\subset S$, and a homomorphism
$\pi_1(S,\xi)\ra\pi_1(\wM_{g,n},f(\xi))$. Then, there exists a projective birational morphism with 
connected fibres $\sg:S'\ra S$, from a normal, irreducible, locally noetherian scheme $S'$, such that:
\begin{enumerate}
\item the morphism $f:U\ra\wM_{g,n}$ extends to a morphism $\ol{f}:S'\ra\wM_{g,n}$;
\item  for all points $P\in S\ssm U$, the curve $C\ra S'$, corresponding to the morphism $\ol{f}$, restricts
over $\sg^{-1}(P)$ to a linear family (see Definition~\ref{linearly rigid}). 
\item the morphism $\sg$ is an isomorphism over $U$ and over neighborhoods of points $P\in S\ssm U$ such
that the set of limit curves of $U$ in $P$ contains a linearly rigid curve.
\end{enumerate}
\end{theorem}
\begin{proof}Proceeding like in the proof of Theorem~\ref{ext}, one constructs a proper birational map
$\sg:S'\ra S$, from a normal, irreducible, locally noetherian scheme $S'$, which is an isomorphism above $U$
and such that the given morphism $f:U\ra\wM_{g,n}$ extends to a morphism $\ol{f}:S'\ra\ccM_{g,n}$. Since
the induced homomorphism  $f_\ast:\pi_1(U,\xi)\ra\pi_1(\wM_{g,n},f(\xi))$ factors through the 
epimorphism $\pi_1(U,\xi)\tura\pi_1(S',\xi)$ induced by the injection $U\hookra S'$, from the remark
following Theorem~\ref{ext grothendieck}, it follows that the domain of $\ol{f}$ is actually contained in
$\wM_{g,n}$. Let us observe that the set of points $P\in S$ such that the set of limit curves of $U$ in $P$
contains a  linearly rigid stable curve is then a Zariski open $U'$ of $S$.

Let $P\in S\ssm U$ and let $C\ra\sg^{-1}(P)$ be the family of $n$-pointed, genus $g$ stable curves of
compact type, corresponding to the restriction of $\ol{f}$ to $\sg^{-1}(P)$. Then, the
induced monodromy representation $\pi_1(\sg^{-1}(P),x)\ra\pi_1(\wM_{g,n},\ol{f}(x))$ is trivial, because
it factors through the monodromy representation $\pi_1(S',x)\ra\pi_1(\wM_{g,n},\ol{f}(x))$ and then
through the homomorphism $\pi_1(S,\xi)\ra\pi_1(\wM_{g,n},f(\xi))$ given in the hypotheses of the
theorem. Thus, by Theorem~\ref{triviality}, the family $C\ra\sg^{-1}(P)$ is linear. In particular, if for
some $x\in\sg^{-1}(P)$, the curve $C_x$ is linearly rigid, the family $C\ra\sg^{-1}(P)$ is trivial and the fibre
$\sg^{-1}(P)$ is contracted by the morphism $\ol{f}:S'\ra\wM_{g,n}$. The same arguments used in the
proof of Theorem~\ref{ext} imply that the morphism $f:U\ra\wM_{g,n}$ extends to a morphism
$f':U'\ra\wM_{g,n}$. By construction, the morphism $\sg:S'\ra S$ has all the required properties.

\end{proof}

\smallskip
\noindent {\bf Acknowledgements}\, This paper grew out of some conversations I had with
Ernesto Mistretta at Jussieu. I thank him for this as well as for other useful remarks he made on 
a preliminary version of the paper. Special thanks go to the referee of IMRN
who pointed out several inaccuracies in a preliminary version of this paper and gave valuable advice.
I thank Pierre Lochak for useful conversations on this subject and
for his steady interest in my mathematical work. Finally, I thank my father whose financial support allowed
me to complete this work.

\bigskip

\noindent Address:\, Institut de Math{\'e}matiques de Jussieu; 175, rue du Chevaleret; 75013 Paris.\\
E--mail:\, boggi@math.jussieu.fr

\end{document}